\documentclass[11pt]{article}
\usepackage{amssymb,amsmath,amsthm,graphics,subfigure}
\makeatletter
\renewcommand{\fnum@figure}[1]{\figurename~\thefigure\ignorespaces}
\makeatother

\newtheorem{theorem}{Theorem}[section]
\newtheorem{proposition}[theorem]{Proposition}
\newtheorem{lemma}[theorem]{Lemma}
\newtheorem{corollary}[theorem]{Corollary}
\newtheorem{conjecture}[theorem]{Conjecture}
\newcommand{\C}{\mathbb{C}}
\newcommand{\Z}{\mathbb{Z}}
\newcommand{\Hom}{\mathrm{Hom}}
\newcommand{\op}{{\mathrm{op}}}
\newcommand{\cop}{{\mathrm{cop}}}
\newcommand{\bd}{\partial}
\newcommand{\cross}{\times}
\newcommand{\tensor}{\otimes}
\newcommand{\D}{{\Delta}}
\newcommand{\eps}{{\epsilon}}
\renewcommand{\t}[1]{{\tilde{#1}}}

\newcommand{\fig}[1]{Figure~\ref{#1}}

\begin{document}

\title{Involutory Hopf algebras and 3-manifold invariants}
\author{Greg Kuperberg \\ UC Berkeley}

\date{May 19, 1990}

\maketitle

\begin{abstract}
We establish a 3-manifold invariant for each finite-dimensional, involutory
Hopf algebra.  If the Hopf algebra is the group algebra of a group $G$, the
invariant counts homomorphisms from the fundamental group of the manifold to
$G$.  The invariant can be viewed as a state model on a Heegaard diagram or a
triangulation of the manifold.  The computation of the invariant involves
tensor products and contractions of the structure tensors of the algebra.  We
show that every formal expression involving these tensors corresponds to a
unique 3-manifold modulo a well-understood equivalence.  This raises the
possibility of an algorithm which can determine whether two given 3-manifolds
are homeomorphic.
\end{abstract}

\section{Introduction}

This paper describes some new invariants for triangulations and
Heegaard diagrams of 3-manifolds.  Specifically, for every
finite-dimensional, involutory Hopf algebra $H$, we define an invariant
$\sharp(M,H)$ for closed, oriented 3-manifolds $M$.  The invariants can be
extended in various ways to invariants of oriented manifolds with
boundary.

We will convert a Heegaard diagram to an expression in terms of the
structure tensors of a Hopf algebra.  The value of the expression is
then the value of the invariant.  Using index notation, we write
$M_{ab}{}^c$ and $\D_a{}^{bc}$ for the multiplication and
comultiplication tensors of some involutory Hopf algebra, and $S_a{}^b$
for the antipode map.  We may consider the set of all expressions in
terms of these tensors, modulo the axioms of a Hopf algebra, which
consist of identities which these tensors satisfy.  We will show that
if two (prime, closed, oriented) 3-manifolds produce equivalent formal
expressions, then they are homeomorphic.  Conversely, we may consider
the question of determining whether two such expressions are formally
equivalent.  For example, is the equation: $$M_{ab}{}^c
S_c{}^d\D_d{}^{ab} = M_{ab}{}^c S_c{}^d \D_d{}^{ba}$$ an identity that
follows directly from the axioms of an involutory Hopf algebra?  The
answer is yes if and only if two certain oriented 3-manifolds are
homeomorphic (in this case $L(3,1)$ and $L(3,2)$).

This work is part of a new area of mathematics which we will
call state model topology, as described in \cite{Kauffman:state,Jones:pacific}.
As the starting point of state model topology, we are
given a topological object, for example a knot, with a combinatorial
description, for example a knot projection.  We wish to find a state
model (or collection of state models) which we can assign to all such
combinatorial descriptions.  We then hope that the partition function
of the state model is a topological invariant, i.e. is independent of
the chosen combinatorial description.  Typically, we demonstrate
invariance under an appropriate set of local moves on the combinatorial
description, for example the three Reidemeister moves.  (Strictly
speaking, the partition function may be a topological ``covariant''
because its value may change in a simple way under some of the moves.)

A \emph{state model} $M$ consists of a (commutative) ring $R$ (usually the
complex numbers), a bipartite graph $G$, the \emph{connectivity graph},
whose vertices are labeled as \emph{atoms} and \emph{interactions}; a set
$S_A$ for each atom $A$, called the \emph{state set} of $A$; and a
function $w_I:A_1\cross A_2 \cross\ldots \cross A_n \to R$ for each
interaction $I$ (where $A_1,\ldots,A_n$ are the neighbors of $I$),
called the \emph{weight function} or the \emph{Boltzmann weights} of
$I$.  A \emph{state} of $M$ is a function $s$ on the atoms of $M$ such
that $s(A) \in S_A$.  The \emph{weight} $w(s)$ of a state $s$ is defined as
the product of the $w_I$'s evaluated at the state $s$ when this product
converges, and in particular when $G$ is finite.  Finally, the
\emph{partition function} $Z(M)$ is defined to be the sum of $w(s)$ over all
states $s$ when this sum converges, and in particular all state sets
are finite.

\begin{figure}[ht]
\begin{center}
\includegraphics{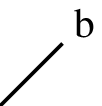}
\caption{\label{f:crossing}}
\end{center}
\end{figure}

To obtain a topological state model, we choose a prescription for
assigning a connectivity graph and weight functions to each
combinatorial description of a topological object.  For example, if we
are given a projection $P$ of a link, we can declare the arcs between
crossings to be atoms and the crossings themselves to be the
interactions.  We choose some $n$-element set $S$ to be the common
state set for all atoms.  If the states of the arcs incident to a given
crossing are labeled $a$,$b$,$c$, and $d$ as in \fig{f:crossing}, we may
define the weight function $w$ of a crossing by:  $$w(a,b,c,d) =
t\delta(a,b)\delta(c,d)+t^{-1}\delta(a,c)\delta(b,d),$$ where $t$ is
chosen so that $n = -(t^2+t^{-2})$, and $\delta(a,b)=1$ when $a=b$ and
0 otherwise.  This state model is a ``link covariant'' called the
Kauffman bracket, which is essentially the Jones polynomial up to
normalization.

Surprisingly, there are many non-trivial topological state models.  We may
start with the general form of such a state model for some kind of
combinatorial description, and interpret combinatorial moves as a list of
constraints on the weights of the model.  Usually there are many more equations
than unknowns, and yet the equations have solutions. For example, the Jones and
HOMFLY polynomials \cite{HOMFLY:invariant} and most of their variants and
generalizations can be defined as topological state models.

We can conclude that the constraints on the weights of a topological state
model have a deep algebraic structure.  Indeed, the constraints imposed by
invariance under the third Reidemeister move are also known as the Yang-Baxter
equations, considered independently in the related area of exactly solved
statistical mechanical models \cite{Baxter:exactly}.  More recently, many
solutions to the Yang-Baxter equations have been discovered via Hopf algebras,
or quantum groups, as is summarized in an inspiring paper by Drinfel'd
\cite{Drinfeld:quantum}.  Recently, Turaev and Reshetikhin have extended this
analysis to produce a 3-manifold invariant which generalizes the Jones
polynomial \cite{RT:manifolds}.

The new invariants are closely related to other topological invariants.
If we restrict to the special case when $M$ is a link complement,
we may convert the invariant $\sharp(M,H)$ to a special case of the machinery
of Drinfel'd, Turaev, and Reshetikhin.  This conversion involves
a bigger Hopf algebra $D(H)$, the quantum double of $H$, which was
defined by Drinfel'd.

The following is an instructive example of a topological state model.
Let $T$ be a triangulation of a connected $n$-manifold $M$.  We orient
the edges of $T$.  We choose a finite group $G$ to be the state set,
and we declare the edges of $T$ to be the atoms.  If we assign $g$ to
an edge with one orientation, we declare this equivalent to assigning
$g^{-1}$ to the same edge with the opposite orientation.  We assign
interactions to the faces (2-simplices) of $T$.   Given a state of the
model, the weight of a face is 1 if the cyclically-ordered product of
the states of the edges (when the edges are oriented in the same
direction around the face) of the face is the identity and 0
otherwise.  We leave it as an exercise that the partition function is
simply $|G|^{v-1} |\Hom(\pi_1(M),G)|$, where $v$ is the number of
vertices of $T$.  Thus, the partition function is an invariant after
dividing by $|G|^{v-1}$.  We will see later that this invariant equals
$\sharp(M,H)$ when $H$ is the group algebra of $G$.

\subsection{Acknowledgments}

Vaughan Jones introduced me to state model topology and Hopf
algebras.   I am also grateful to my advisor, Andrew Casson; Rob Kirby;
and Oleg Viro for their encouragement and their knowledge of 3-manifold
topology.

\section{Tensor notation}

Let $V$ be a finite-dimensional vector space over a field $k$.  For most of
this paper, a tensor will mean an element of some tensor product space $V_1
\tensor V_2 \tensor \ldots \tensor V_n$, where each $V_i$ is either $V$ or
$V^*$.  We will adopt index notation for tensors as it is defined in
\cite{Wald:general}.  In index notation, a tensor $T \in V_1 \tensor \ldots
\tensor V_n$ is written as $T_a{}^b{}^{\ldots}_{c\ldots}$, with $n$ indices
total, where the $k$th index is a superscript (or contravariant index) if $V_k
= V$ and a subscript (or covariant index) if $V_k = V^*$.  The indices should
be distinct letters which have no meaning independent of their use as a
description of the tensor $T$. For example, a vector is written $v^a$, a dual
vector is $w_a$, a linear operator is $L_a{}^b$, a bilinear form is $g_{ab}$,
and so on. If a tensor $T$ has $k$ superscripts and $n-k$ subscript, we say
that $T$ has \emph{type} $(k,n-k)$, and \emph{total type} $n$.

The indices of a tensor can be used for three different operations on
tensors:

\begin{description}
\item[1.] We denote the factor-switching map $a
\tensor b \mapsto b \tensor a$ and its cousins for higher tensor
products by permuting the corresponding indices.  For example, if
$g_{ab} = g_{ba}$, then $g$ is a symmetric bilinear form.

\item[2.] We write the tensor product of two tensors by juxtaposing the
two tensors and giving them distinct indices.  For example, one may
write $L_a{}^b = w_av^b$ if $L$ is a linear operator of rank 1.

\item[3.] The canonical map from $V^* \tensor V$ to $k$ which is called
the trace map, the contraction map, or the evaluation map is denoted by
a repeated index.  The index must appear exactly once as a superscript
and once as a subscript.  Used together, contractions and tensor
products subsume most of the usual operations in linear algebra.  For
example, $w_av^a$ is the value of $w$ at $v$, $L_a{}^bv^a$ is the
operator $L$ applied to $v$, $L_a{}^a$ is the trace of $L$, and so on.
\end{description}

Although index notation only allows basis-independent
computations with tensors, the notation is motivated by computations in
a specific basis.  Suppose the vector space $V$ has dimension $n$ and
has a given basis.  We may interpret the notation $v^a$ to mean the
coordinates of the vector $v$, with $a$ an integer from 1 to $n$, and
in general we may interpret $T_a{}^b{}^{\ldots}_{c\ldots}$ to be the
matrix entries of the tensor $T$.  In this case a tensor formula in
index notation may be interpreted as a formula for the matrix entries
of the tensors, and the two formulas are equivalent given the
convention that we sum over any repeated indices (the Einstein
summation convention).

We may use the summation convention to arrive at a re-interpretation of
certain state models.  Suppose that $M$ is a state model such that
every atom is adjacent to precisely two interactions, and suppose for
simplicity that there is only one state set $S$.  We may declare $S$ to
be the index set of a basis (or perhaps $S$ is the basis) of a vector
space $V$.  We declare each atom $a$ to be an index letter.  We may
then re-interpret the values of each weight function $w_I$ as the
matrix entries of a tensor $T_I$ of total type $n$.  We arbitrarily
declare each of the indices of each $T_I$ to be covariant or
contravariant, provided that each atom is used once as a superscript
and once as a subscript.  In this case the partition function $Z(M)$ is
the result of tensoring together all of the interaction tensors and
then contracting all of the (atom) indices in the natural manner.

Thus, we define a \emph{tensorial state model} to be a list of tensors
over one vector space, together with a formula which is some product of
these tensors with all indices contracted.  The partition function is
defined to be the value of this formula.  More generally, we say that
$M$ is a \emph{tensorial state model with boundary} if $M$ is a formula
which is a product of tensors with some subset of the indices
contracted.  We call the uncontracted indices the boundary of the
model, and the contracted indices the interior.  We define the
\emph{partition tensor} $Z^{a\ldots}_{b\ldots}(M)$ to be the value of the
formula.  The partition tensor $M$ should not be confused with the
partition function of $M$.

The state models defined in this paper will all be tensorial state
models.

To conclude this section, we present a more visual version of index
notation which we will actually use.  To avoid using many different
letters for contracted indices, we may draw each tensor with an inward
arrow for every subscript and an outward arrow for every superscript.  For
example, we write the tensor $T_{ab}{}^{cd}$ as follows:
$$\includegraphics{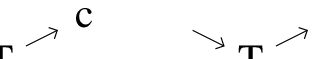}$$
We denote the contraction of two indices by connecting the corresponding
arrows.  For example, $M_{ab}{}^cv^aw^b$ and $M_{ab}{}^b$ are written as:
$$\includegraphics{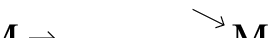}$$
For consistency, we should interpret the strange-looking expression:
$$\includegraphics{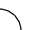}$$
as the scalar $\dim V$.

Arrows can and will cross.  Strictly speaking, the visual notation can
be ambiguous.  For example, it is not clear if the following denotes
$T_{abc}$ or $T_{cab}$, or perhaps even $T_{acb}$:
$$\includegraphics{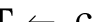}$$
We declare that indices should be read off in counter-clockwise order.  The
meaning of a diagram is still ambiguous if the arrows incident to some tensor
are symmetric under rotation.  However, all such tensors in this paper will be
symmetric under cyclic permutation of their indices, so this ambiguity will be
irrelevant to us.

\section{\label{Hopf}Hopf algebras}

A finite-dimensional \emph{Hopf algebra} is a finite-dimensional vector
space $H$ (over some ground field $k$) such that $H$ and $H^*$ are both
algebras with unit, and the two algebra structures satisfy certain
compatibility axioms.  We will describe the axioms with the fancy
notation of the previous section.  Firstly, we may describe the algebra
structure on $H$ by a \emph{multiplication tensor} $M$ defined so that
$$\includegraphics{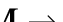}$$
is the product of $v$ and $w$.  $M$ satisfies the associativity axiom:
$$\includegraphics{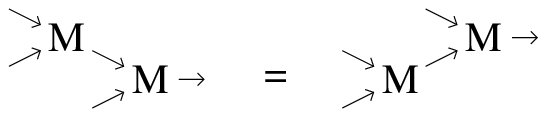}$$
Similarly, the algebra structure on $H^*$ is described by the
\emph{comultiplication tensor} $\D$ which is coassociative:
$$\includegraphics{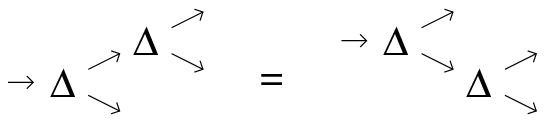}$$
In addition, there are three other tensors $i$, $\eps$, and $S$
called the \emph{unit}, the \emph{counit}, and the \emph{antipode} of
$H$.  The five tensors together satisfy these axioms:
$$\includegraphics{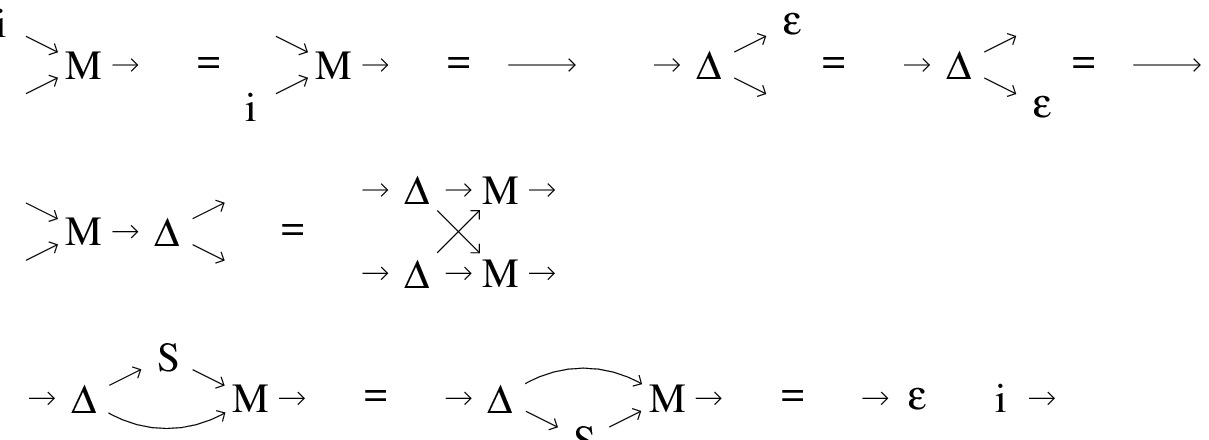}$$
We also assume as an axiom that $S$ has an inverse, although this fact follows
from the other axioms in the finite-dimensional case \cite{Sweedler:hopf}.
Given $M$ and $\D$, the other tensors are unique when they exist, so we may say
that $M$ and $\D$ are the ``meat'' of a Hopf algebra.  A \emph{morphism of Hopf
algebras} is defined to be a linear transformation of the underlying vector
spaces which commutes with the five structure tensors.

Since we have not defined infinite-dimensional Hopf algebras and will
not use them, we will assume that a Hopf algebra is finite-dimensional
unless explicitly stated otherwise.

We adopt the following abbreviations in light of
the associativity and coassociativity axioms: 
$$\includegraphics{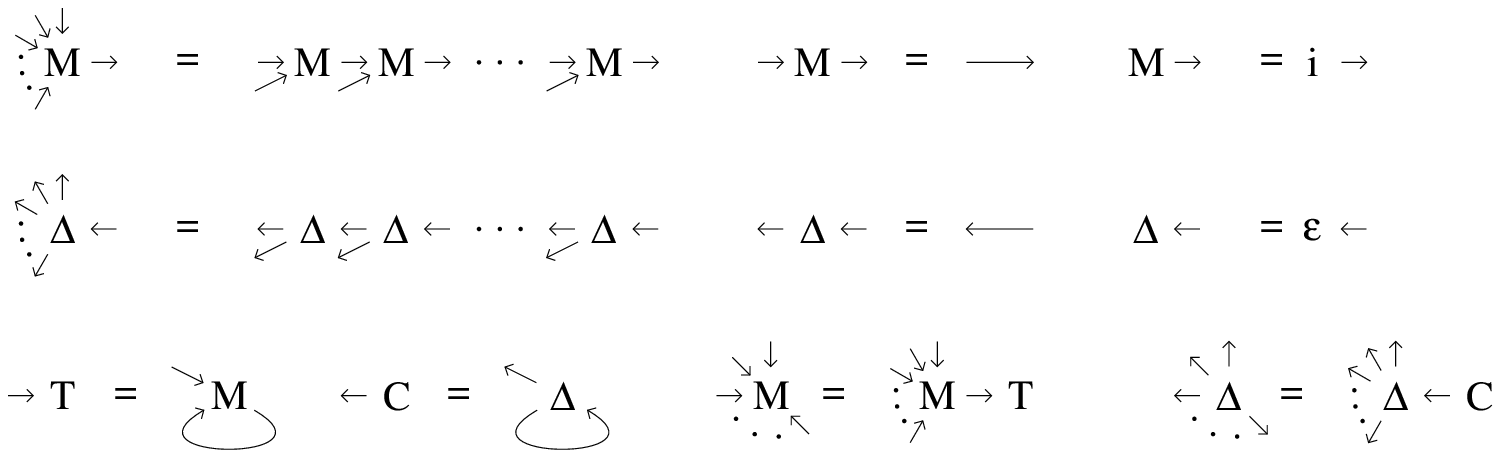}$$
The dual vector $T$ is called the \emph{trace}, and the vector $C$ is
the \emph{cotrace}.  We will call the last two tensors the \emph{tracial
product} and \emph{tracial coproduct}.  Like the diagrams, the tensors
themselves are cyclically symmetric.

One motivation for the definition of a Hopf algebra is the
observation that every group algebra is naturally a Hopf algebra.  Let $G$  
be a finite group.  We may let the elements of $G$ be both
a basis and a dual basis for a vector space $k[G]$, and define
multiplication and comultiplication by:
$$\includegraphics{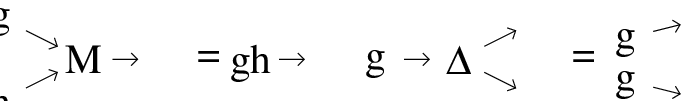}$$
The antipode is induced by the group inverse.  This construction has the nice
property that the morphisms between the Hopf algebras of two groups correspond
to the homomorphisms between the groups themselves.

In general, the antipode plays the role of the group inverse in a Hopf
algebra.  However, the antipode need not satisfy $S^2 = I$. A Hopf
algebra is called involutory if this equation is satisfied.  Many
important Hopf algebras are not involutory, and in particular those
used in the analysis of the Jones polynomial are non-involutory
\cite{Drinfeld:quantum}.  However, we will only obtain results about involutory
Hopf algebras in this paper.

Finite-dimensional Hopf algebras possess an eight-fold duality.  If $H$
is a Hopf algebra, we obtain a new Hopf algebra $H^\op$ on the same
vector space by reversing multiplication:
$$\includegraphics{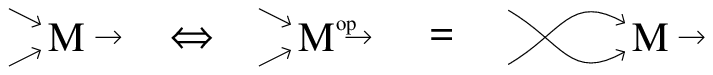}$$
We will show below
that $S^{-1}$ is the antipode of $H^\op$.  Similarly, we define
$H^\cop$ by reversing comultiplication.  Finally, we define the
\emph{dual} of $H$, a Hopf algebra structure on the vector space $H^*$, by
switching $M$ with $\D$ and $\eps$ with $i$ and reversing all the
arrows (Actually, what we are calling $H^*$ is usually defined as
$H^{*,\op,\cop}$, and vice-versa):  $$\includegraphics{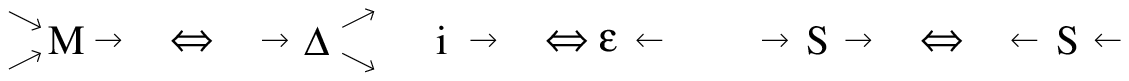}$$

We establish some basic identities involving Hopf algebras.  The proof
of each identity will use the following observation:

\begin{lemma} The tensors:
$$\includegraphics{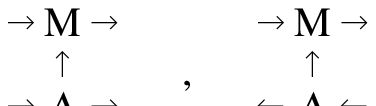}$$
when viewed as vector space endomorphisms of $H \tensor H$ and
$H \tensor H^*$, are invertible. \label{ladder}
\end{lemma}
\begin{proof}
We compute:
$$\includegraphics{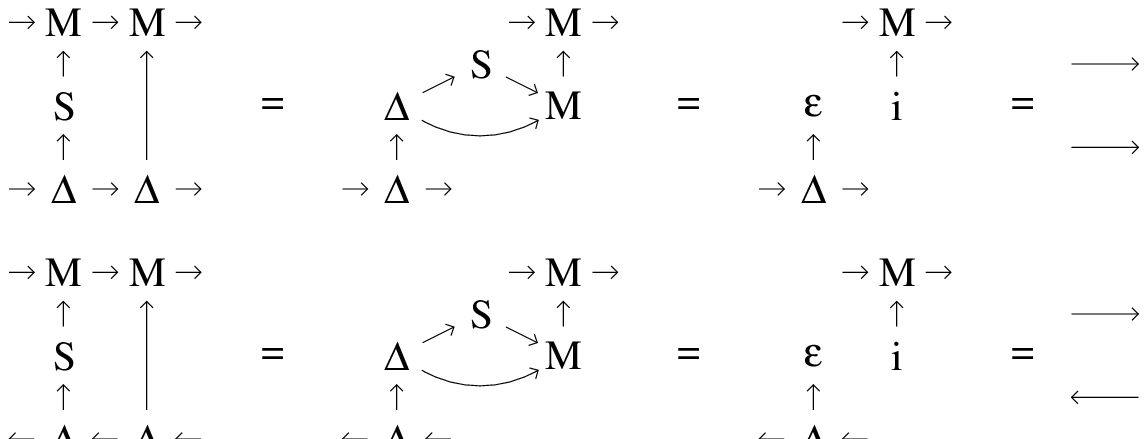}$$
\end{proof}

\begin{lemma} The following identities hold in any Hopf
algebra:
$$\includegraphics{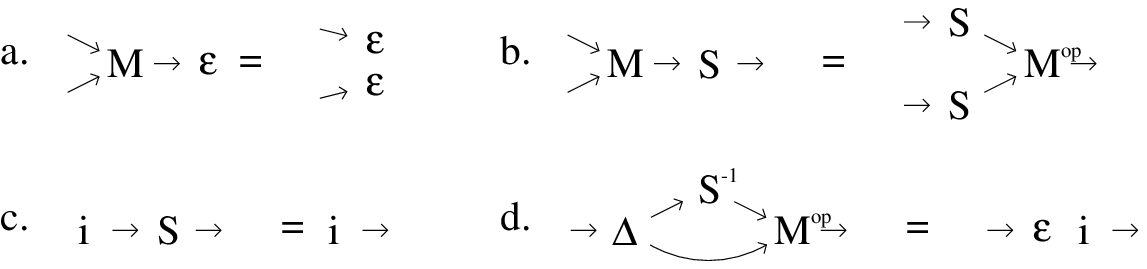}$$ \label{basic}
\end{lemma}
\begin{proof}
We compute: $$\includegraphics{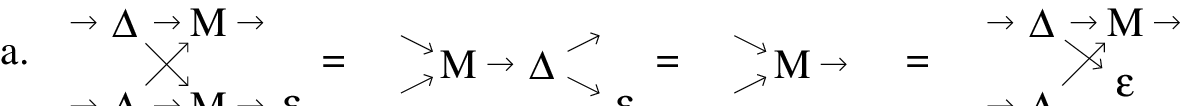}$$
\end{proof}

We will need a few other definitions involving Hopf algebras.  If $H$ is a Hopf
algebra, a \emph{left integral} of a Hopf algebra is a dual vector $\mu_L$ that
satisfies: $$\includegraphics{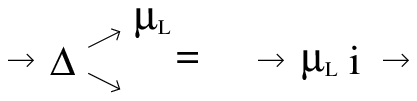}$$ One can similarly define a
\emph{right integral}, a \emph{left cointegral}, and \emph{right cointegral}:
$$\includegraphics{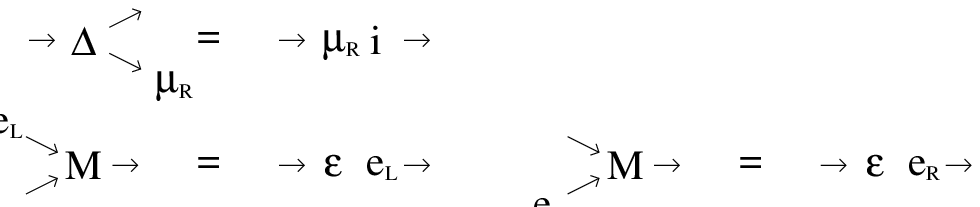}$$
A Hopf algebra always has a non-zero integral and cointegral.  By
Lemma~\ref{basic}, $\eps$ is an algebra homomorphism. In general, if $A$ is a
finite-dimensional algebra over $k$ and $\eps:A \to k$ is any homomorphism,
then there exists a non-zero element $e$ such that $ae = \eps(a)e$ for all $a
\in A$ \cite{Bergman:personal}.

A Hopf algebra $H$ is called \emph{semisimple} if it is semisimple as an
algebra, and cosemisimple if $H^*$ is semisimple.  Equivalently,
$H$ is semisimple if the tensor: $$\includegraphics{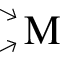}$$
is a non-degenerate bilinear form.  We say that $H$ has \emph{invertible
dimension} if $\dim H \ne 0$ in the ground field of $H$.

In an involutory Hopf algebra, the trace $T$ is also a left and right integral
\cite{LR:semisimple0}.  To prove this, we use the fact that there exists some
non-zero integral, and we show that $T$ is a multiple of it:

\begin{lemma}[Radford and Larson] The tensor:  $$\includegraphics{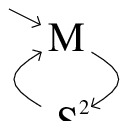}$$
is always a right integral (which may be zero). \label{tracelemma}
\end{lemma}
\begin{proof}
Let $\mu_R$ be a non-zero right integral and let $e_L$ be a non-zero
left cointegral.  We first observe that: $$\includegraphics{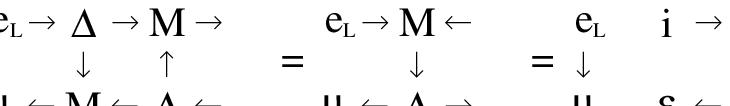}$$
Therefore:
$$\includegraphics{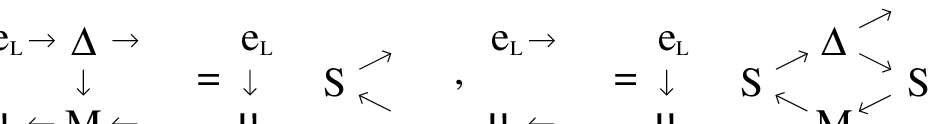}$$
by Lemma~\ref{ladder}.  Since $\mu_R$ and $e_L$ are non-zero, the
last equation shows that $\mu_R(e_L)$ must be non-zero.  Finally,
we obtain:
$$\includegraphics{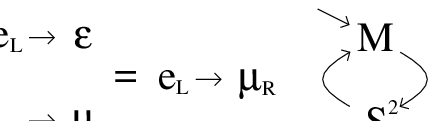}$$
as desired.
\end{proof}

By virtue of the fact that: $$\includegraphics{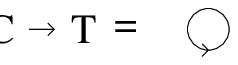}$$
we obtain:
\begin{corollary} If $H$ is involutory, then:
$$\includegraphics{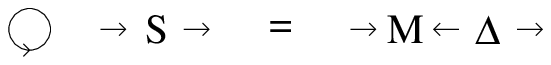}$$
\label{tracegen} \end{corollary}
Thus, Lemma~\ref{tracelemma} demonstrates not only that $T$ is
a right integral, but also that $H$ is semisimple and cosemisimple.
Moreover, if $H$ has invertible dimension, we can express all of the
structure tensors of a Hopf algebra in terms of tracial product and
coproduct!

\section{Heegaard diagrams and triangulations}

In the section (and the rest of the paper) we will work in the smooth
category of manifolds.  We will use the term $n$-manifold to mean
an oriented (not merely orientable) and compact $n$-manifold which may
have boundary and may or may not be connected.

We recall some basic elements of 3-manifold topology as described in
\cite{Hempel:book}.  A \emph{Heegaard splitting} is a decomposition of a closed
3-manifold into two handlebodies of the same genus glued along their
boundary; the boundary is called a \emph{Heegaard surface}. Given such a
decomposition, we may arbitrarily label the two handlebodies upper and
lower.  If the handlebodies are labeled, then an orientation of the
Heegaard surface induces an orientation of the manifold:  At a point
$p$, we complete a positively-oriented basis $(e_1,e_2)$ of the tangent
space of $S$ at $p$ to a (positive) basis $(e_1,e_2,e_3)$ of $TM|_p$ by
choosing $e_3$ to point from the upper handlebody to the lower one.

If $S$ is a surface, we define a \emph{handlebody diagram} $d$ on $S$
to be a prescription for gluing a handlebody to $S$.  Specifically, it
is a family of disjoint circles $\{c_i\}$ which divides $S$ into planar
regions.  (Often it is assumed that $\{c_i\}$ does not separate $S$,
but we do not want this hypothesis.)  We obtain a handlebody $H_d$ from
$d$ by the following procedure:  We start with $S \cross I$ and glue a
2-handle along a tubular neighborhood of each circle $c_i \cross
\{0\}$.  The boundary of the result is $S \cross \{1\}$ plus a number of
spheres.  We eliminate the spherical boundary components by gluing in
balls to obtain $H_d$.

We say that $D$ is a \emph{Heegaard diagram} of a closed 3-manifold on a
surface $S$ is a prescription for a Heegaard splitting of a manifold
$M_D$, i.e. it is a pair of handlebody diagrams $d_l$ and $d_u$ which
describe the upper and the lower handlebodies of the splitting.  We
assume that $d_l$ and $d_u$ are transverse.

We will need to generalize the definition of a Heegaard diagram to
account for manifolds with boundary, and to relax the constraint that
the circles of a handlebody diagram divide the surface of the diagram
into planar regions.  To this end, we also label the boundary of a
3-manifold $M$ by arbitrarily dividing the boundary components of $M$
into two disjoint subsets:  The upper boundary and the lower boundary.
We also define a \emph{puncture move} on a 3-manifold $M$ with labeled
boundary to consist of removing a ball from $M$ and labeling the new
boundary as upper or lower.  We define a \emph{bordism} to be an
equivalence class of 3-manifolds with labeled boundary under the
puncture move.  By abuse of terminology, we will not distinguish
between a bordism and a representative of a bordism.

A \emph{compression body} is a manifold which is obtained from $S \cross
I$ by gluing 2-handles on one side, where $S$ is a surface (usually one
requires that all spherical boundary components of a compression body
are capped by balls, but in the context of bordisms this does not
matter).  A compression diagram is then a generalization of a
handlebody diagram:  It is a prescription for gluing a compression body
to a surface $S$ consisting of a collection of disjoint circles.
However, the circles are not required to divide the surface into planar
regions.  We define a generalized Heegaard diagram to be a pair of
compression diagrams (an upper diagram and a lower diagram) on one
surface.

If we are given a generalized Heegaard diagram $D$ on $S$, we may form a
bordism $M$ by gluing together the corresponding compression
bodies along $S$.  The boundary of the lower compression body becomes 
the lower boundary of the bordism, and the same for the upper
compression body.

We will need a set of moves to convert any Heegaard diagram of a
bordism to any other diagram for the same bordism.  Also, it may not be
obvious that every bordism has a Heegaard diagram.  We may solve both
of these problems using the theory of Morse functions, as described in
\cite{Cerf:strat}.  We summarize the relevant parts of the theory:  If $M$ is a
bordism, we first choose a Riemannian metric for $M$.  We choose a
Morse function $f$ on $M$ which attains its minimum on the lower
boundary of $M$ and its maximum on the upper boundary of $M$.  Using
the Riemannian metric, we can consider the flow $\phi$ on $M$ given by
the gradient of $f$.  For each fixed point $p$ of $\phi$, we consider
the descending manifold of $p$, the set of all points that flow into
$p$.  The descending manifold is necessarily a disk of some dimension,
and in this way we obtain a handle decomposition of $M$.  If $f$ is in
general position, each $n$-handles is attached only to $m$-handles with
$m<n$ (not merely $m \le n$) and possibly to the lower boundary.
Finally, it is easy to convert such a handle decomposition to a
generalized Heegaard diagram:  We let the lower compression body be the
union of the 0-handles, the 1-handles, and a tubular neighborhood of
the lower boundary, and we let the upper compression body be the
complement of the lower compression body.  We obtain an upper circle
for each 2-handle by intersecting the core (descending manifold) of the
2-handle with the Heegaard surface, and similarly we intersect the
cocore (ascending manifold) of each 1-handle with the Heegaard surface
to obtain the lower circles.

Since every bordism has a Morse function, every bordism has
a Heegaard diagram.  On the other hand, given a Heegaard diagram,
it is easy to construct a Morse function which reproduces it.
We can therefore obtain a set of moves on Heegaard diagrams from
a set of moves on Morse functions.  The upshot is the following theorem:

\begin{figure}[ht]
\subfigure[]{\label{f:isomove}\includegraphics{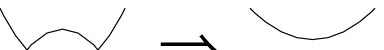}}
\subfigure[]{\label{f:handlemove}\includegraphics{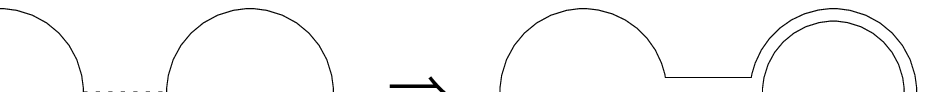}}
\subfigure[]{\label{f:stablemove}\includegraphics{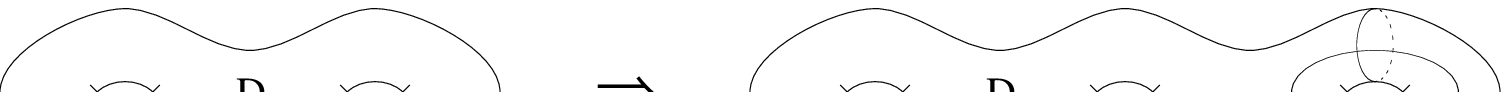}}
\caption{ }
\end{figure}

\begin{theorem} Let $D$ be a generalized Heegaard diagram on a
surface $S$ and consider the following moves on $D$:
\begin{description}
\item[1.]  Homeomorphism of the diagram.  Using a homeomorphism from a surface
$S$ to a surface $T$, we carry $D$ to a diagram on $T$.
\item[2.]  The two-point move.  We isotop the lower circles of $D$ relative
to the upper circles.  If this isotopy is in general position, it reduces
to a sequence of two-point moves, as shown in \fig{f:isomove}.
\item[3.]  Sliding one circle past another.  Suppose $C_1$ and $C_2$ are
two circles of $D$, both lower or both upper, and $A$ is an arc
which connects $C_1$ to $C_2$ but does not cross any other circle.
We add a detour to $C_1$ that follows $A$, goes around $C_2$,
and then goes back along $A$, as shown in \fig{f:handlemove}.
\item[4.]  Creating a trivial circle.  We add a contractible
circle to $D$ which is disjoint from all other circles of $D$.
\item[5.]  Stabilization.  We remove a disk from $S$ which is disjoint
from all circles of $D$ and replace it by a punctured torus with
one lower circle and one upper circle, as shown in \fig{f:stablemove}.
\end{description}
\end{theorem} \begin{proof} Given a bordism $M$, we can construct a
Heegaard diagram from a function $f$ if $f$ has three properties:  $f$
is a Morse function, the descending manifold of each critical point of
$f$ avoids other critical points of the same degree, and the descending
manifolds index 2 critical points are transverse to the ascending
manifolds of the index 1 critical points.  Following \cite{Cerf:strat}, the
space ${\cal F}$ of all smooth functions on $M$ has a stratification
${\cal F}_0, {\cal F}_1, {\cal F}_2, \ldots$ in which ${\cal F}_0$ is
the set of functions $f$ with all three properties.  ${\cal F}_1$ is a
codimension 1 subspace of ${\cal F}$ which is divided into three
pieces:  In each piece, an element $f$ fails to have one of the three
desired properties.  Given two elements $f_0$ and $f_1$ of ${\cal F}_0$,
we choose a path of functions $f_t$ in general position.  The
path avoids the codimension 2 subspace ${\cal F}_2$, and is transverse
to ${\cal F}_1$.  Each time it crosses ${\cal F}_1$, the effect on the
corresponding Heegaard diagram is one of the above moves.  For example,
if $f_t$ fails to be a Morse function, then a pair of critical points
is created or destroyed, the Heegaard diagrams corresponding to
$f_{t+\eps}$ and $f_{t-\eps}$ differ by stabilization or creation of a
trivial circle.  If a descending manifold of a critical point lands on
another critical point of the same degree, the corresponding move is
sliding one circle past another.  Finally, if a descending manifold and
an ascending manifold fail to be transverse, the effect is a two-point
move. \end{proof}

In order to say that our state models can be defined on a triangulation
of a 3-manifold, we construct a canonical way to convert a triangulation
to a Heegaard diagram.  Since a triangulation is a kind of handle
decomposition, we can let the lower handlebody be a tubular
neighborhood of the 1-skeleton as usual.  We add a lower circle for
each edge and an upper circle for each face.  An intersection between
a lower circle and an upper circle corresponds to an ordered pair
consisting of an edge and a face that contains the edge.

\section{\label{invariant}The invariant and the proof that it works}

Let $D$ be a (generalized) Heegaard diagram on a surface $S$.  We assume
that $D$ is an oriented diagram, i.e. its circles are oriented (we
include reversing the orientation of a circle as a move on oriented
Heegaard diagrams).

Let $H$ be a finite-dimensional involutory Hopf algebra of invertible
dimension.  We define the quantity $\sharp(D,H)$ as follows:  To each upper
circle $l$, we assign the tensor:
$$\includegraphics{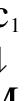}$$
where the indices
$c_1,\ldots,c_n$ correspond to the crossings on $l$ in the order that
they are encountered if we travel along $l$ in the positively oriented
direction.  We assign the tensor:
$$\includegraphics{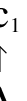}$$
to each lower circle in the
same way.  If, at a given crossing $c$, the tangent vectors of the
lower circle and the upper circle, in that order, form a positively
oriented basis for $TS$ at $c$, we contract the tensors assigned to the
circle at the index corresponding to $c$.  If the vectors form a
negatively oriented basis, we interpose the antipode map before
contracting:
$$\includegraphics{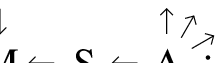}$$
Let $Z(H)$ be the result of this series of tensor
products and contractions. (Note that if we choose a basis for $H$, the
summations corresponding to the contractions can alternatively be
viewed as a sum over states, with $Z(H)$ the resulting partition
function.) Define $$\sharp(D,H) = Z(H) (dim H)^{g(S) - n_u - n_l},$$ where
$n_u$ is the number of upper circles, $n_l$ is the number of lower
circles, and $g(S)$ is the genus of $S$.

\begin{theorem} $\sharp(D,H)$ is an invariant of bordisms.
In particular, it is an invariant of closed 3-manifolds.
\end{theorem}
\begin{proof}  We show invariance under the different moves:

1) Adding a trivial circle.  This has the effect of contributing
the following factor to $Z(H)$: $$\includegraphics{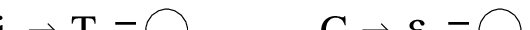}$$
The factor is canceled by the normalization of $\sharp(D,H)$.

2) Stabilization.  This contributes the following factor to $Z(H)$:
$$\includegraphics{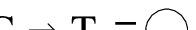}$$
which is also canceled by normalization.

3) Orientation reversal.  An upper circle corresponds to a factor of
$Z(H)$ of the form:
$$\includegraphics{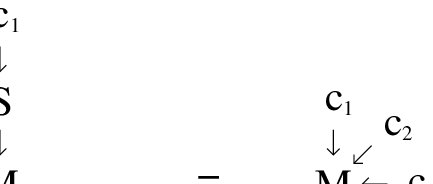}$$
The equality, which follows from Lemma~\ref{basic}, demonstrates invariance. 
The same argument works for lower circles.

4)  The two-point move. If the two circles are oriented properly, part
of the expression for $Z(H)$ is of the form: $$\includegraphics{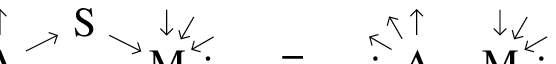}$$
The equality demonstrates invariance.

\begin{figure}[ht]
\begin{center}
\includegraphics{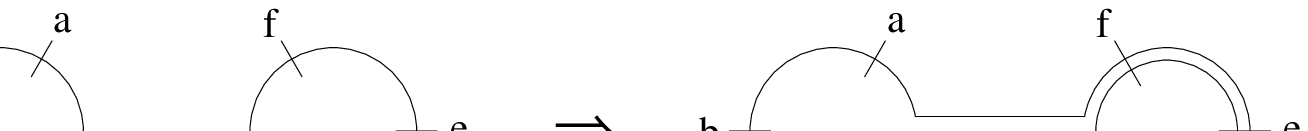}
\caption{\label{f:slide}}
\end{center}
\end{figure}

5)  Sliding a circle.   We can assume that we are sliding one lower
circle past another lower circle (the corresponding move for
upper circles is equivalent by duality).  We assume as a representative
case that both circles have three crossings.  \fig{f:slide} depicts the
move after a possible change of orientation of the circles involved.
The ``before'' picture corresponds to the following factor of $Z(H)$:
$$\includegraphics{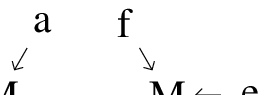}$$
and the ``after'' picture corresponds to the following replacement for this
factor:
$$\includegraphics{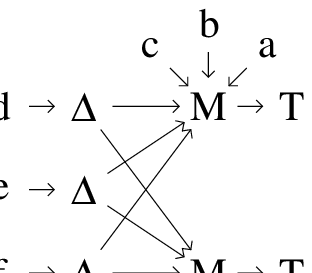}$$
The following algebra demonstrates the equality of these two expressions:
$$\includegraphics{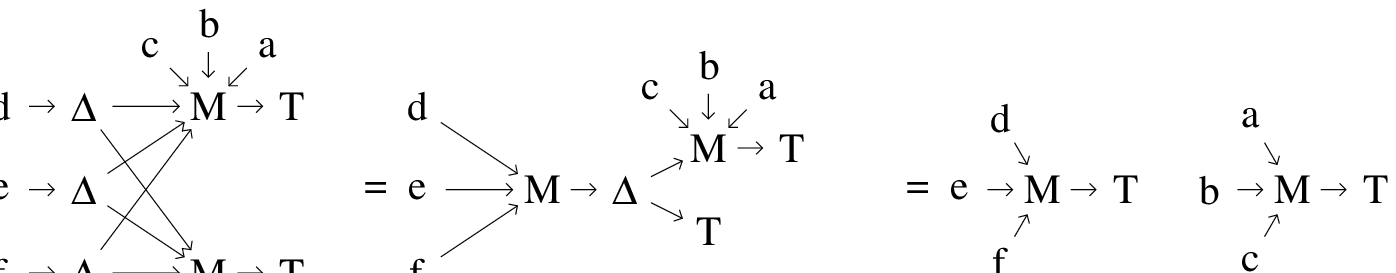}$$
This completes the proof of invariance. 
\end{proof}

The theorem warrants the notation $\sharp(M,H) = \sharp(D,H)$, where
$M$ is a bordism and $D$ is a generalized Heegaard diagram for $M$.

\section{\label{prop} Properties of the invariant and special cases}

We first observe that for a closed manifold $M$, $\sharp(M,H) =
\sharp(M,H^*)$, and $\sharp(-M,H) = \sharp(M,H^\op) =
\sharp(M,H^\cop)$, where $-M$ is $M$ with the oppose orientation.
Also, for any two bordisms, $M_1$, and $M_2$, $\sharp(M_1 \sharp M_2,H)
= \sharp(M_1,H)\sharp(M_2,H)$, because we can choose a diagram for $M_1
\sharp M_2$, the connected sum of $M_1$ and $M_2$, which is a connected
sum of diagrams of $M_1$ and $M_2$.

We will investigate the invariant $\sharp(M,H)$ in the special case
that $H$ is a group algebra, and in the special case that $M$ is a link
complement.

Let $G$ be a finite group and let $k[G]$ be its group algebra.  As
before, view the elements of $G$ as a basis for $k[G]$.  In this basis,
the matrix entry $\D^{ab}_c$ is 1 when $a = b = c$ and 0 otherwise,
while $M_{ab}^c$ is 1 when $ab = c$ (as group elements), and 0
otherwise.  The trace $T^a$ equals $\dim H_G = |G|$ when $a = e$, the 
identity, and 0 otherwise.  The cotrace $C_a$ is 1 for all $a$.  Thus,
the tracial product $M_{a_1a_2\ldots a_n}$ equals $|G|$ when $\prod a_i
= e$ and 0 otherwise, while the tracial coproduct $\D^{a_1a_2\ldots
a_n}$ is 1 when all indices are equal and zero otherwise.

Let $M$ be a closed manifold or a manifold with some upper boundary
(but no lower boundary).  Viewing the invariant $\sharp(M,k[G])$ as a
state model on a Heegaard diagram, we observe that the weight of a
state is zero unless for each lower circle the states assigned to all
crossings on the circle are equal.  Thus, we may say that the states
are assigned to the lower circles themselves.  If the Heegaard diagram
comes from a triangulation of $M$, the lower circles correspond to
edges, so we may say that states are assigned to the edges.  Meanwhile,
the interaction on the upper circles dictates that a state has zero
weight unless the product of the group elements around a face is the
identity.  Thus, the state model reduces to the one described in the
introduction.

\begin{figure}[ht]
\begin{center}
\scalebox{.2}{\includegraphics{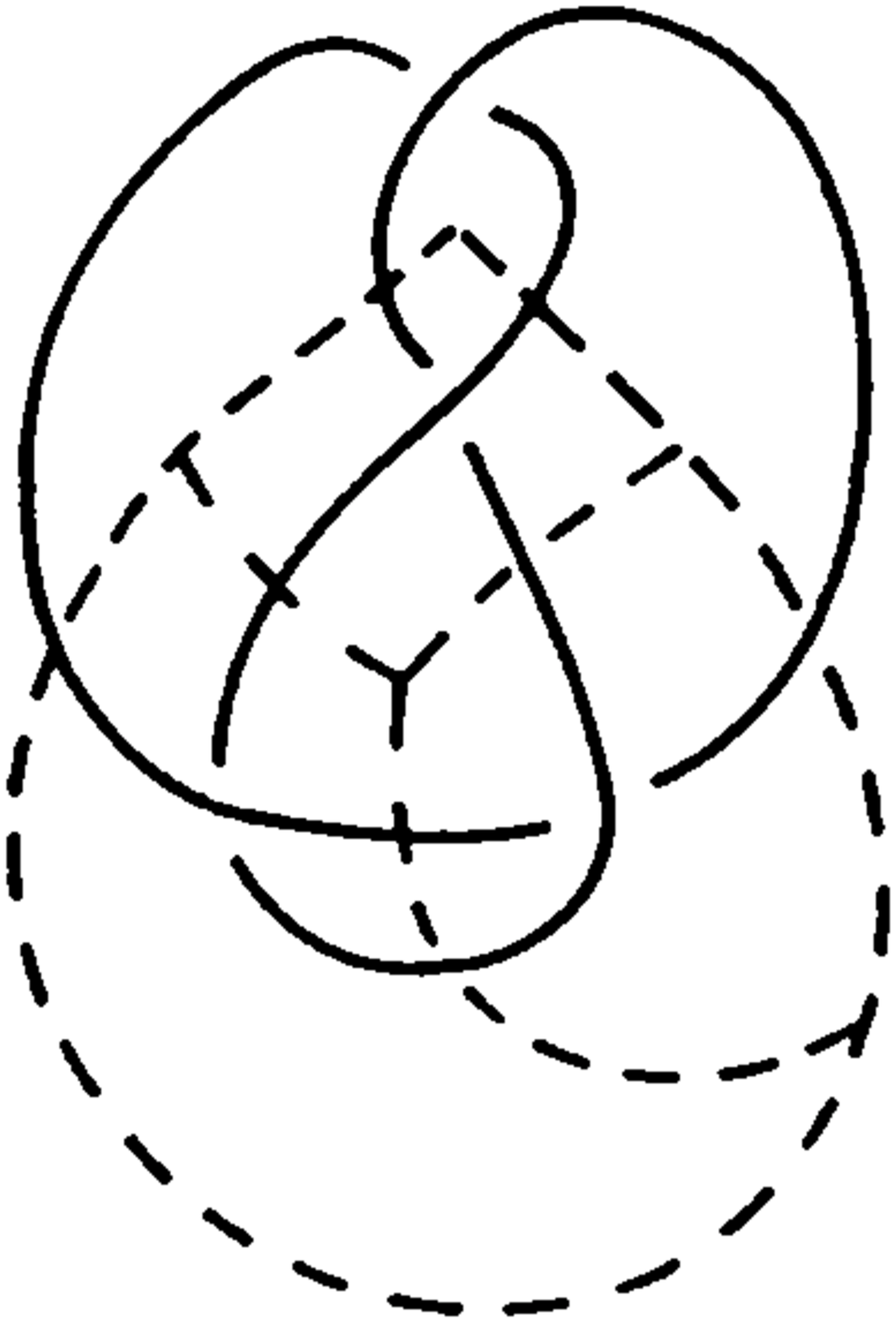}}
\caption{\label{f:squares}}
\end{center}
\end{figure}

\begin{figure}[ht]
\begin{center}
\subfigure[]{\includegraphics{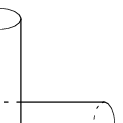}}
\subfigure[]{\includegraphics{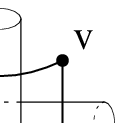}}
\subfigure[]{\includegraphics{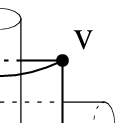}}
\subfigure[]{\includegraphics{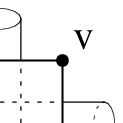}}
\caption{\label{f:square}}
\end{center}
\end{figure}

Let $M_L$ be the complement of (an open neighborhood of) a link $L$,
where we declare that the boundary of $M_L$ is upper boundary.  We wish
to reduce $\sharp(M_L,H)$ to a state model on a projection $P$ of $L$
(we assume $P$ is a projection of $L$ onto the sphere rather than the
plane).  To construct a convenient handle decomposition for $M_L$, we
first consider the tiling of the sphere by squares which is dual to the
link projection, as illustrated in \fig{f:squares}.  There will be one
3-handle $u$, and one 0-handle $v$; imagine $u$ as the region above $P$
and $v$ as a big, fat vertex below $P$.  \fig{f:square}(a) shows a piece of
$M_L$ corresponding to a given square in \fig{f:squares}.  We add a 1-handle
for every edge in the tiling by squares, all attached to the lone
0-handle $l$, as in \fig{f:square}(b).  Next, for each square we add a
2-handle which is attached to the 1-handles on the left and right sides
in \fig{f:square}(c).  Finally, we add a second 2-handle which is attached to
all four neighboring 1-handles, as in \fig{f:square}(d).  As before, we
convert the handle decomposition to a Heegaard diagram by taking the
lower handlebody to be a tubular neighborhood of the 1-skeleton.

Our goal is a model in which states are assigned to the arcs between
crossings of $P$ and interactions are assigned to the crossings
themselves.  Equivalently, we can orient the link $L$ and assign a
tensor of type (2,2) to each crossing, and wherever two crossings are
connected, we contract the corresponding tensors.  We convert the
assembly of handles in \fig{f:square} to tracial product and coproduct
tensors and we modify the coproduct tensors to obtain a tensor of type
(2,2), to obtain the following replacement rule for a left-handed
crossing: 
$$\includegraphics{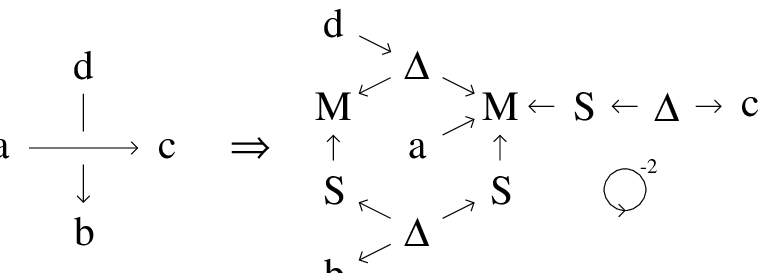}$$
Using Corollary~\ref{tracegen}, we may simplify this tensor:
$$\includegraphics{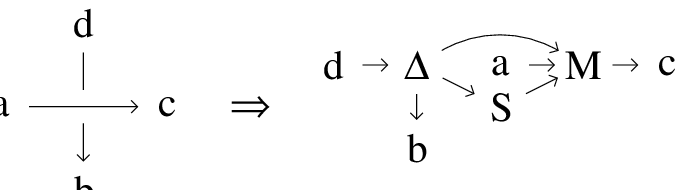}$$

In \cite{KR:links}, Kirillov and Reshetikhin describe a link invariant for
any representation of any Hopf algebra with a special property called
quasi-triangularity.  In particular, for any Hopf algebra $H$
(involutory or not), there is a bigger Hopf algebra $D(H)$, called the
quantum double of $H$, which can be defined on the vector space $H
\otimes H^*$ and which is always quasi-triangular.  The algebras $H$
and $H^*$ are contained in $D(H)$ by the inclusions $a \to a \otimes
\eps$ and $a \to a \otimes i$.  Consider the representation of $D(H)$
on itself given by left multiplication, and let $\rho$ be the
restriction of this representation to $H$ by the inclusion $a \to a
\otimes T$, where the trace $T$ is interpreted as an element of $H^*$.
In this context it is easy to show that Reshetikhin and Turaev's
invariant in the special case of the representation $\rho$ of $D(H)$ is
the same as our invariant $\sharp(M,H)$.  Presumably our invariant, in
its full generality, is a special case of Reshetikhin and Turaev's
extension of this invariant to 3-manifolds.

\section{Tensor notation revisited}

To make the invariant $\sharp(M,H)$ look good, we define a formal calculus of
expressions involving structure tensors of a Hopf algebra, or more
generally arbitrary expressions involving tensor products and
contractions of tensors.  We start by defining a contraction
category as a natural setting for such computations.

A \emph{contraction category} is a set $T$, whose elements may be called
tensors, with additional structure consisting of two functions $i$ and
$o$ from $T$ to $\Z^+$, where $i(t)$ is called the \emph{in-degree} of
$t$ and $o(t)$ is called the \emph{outdegree} (we will also say that $a$
has \emph{type} $(o(t),i(t))$), and a function $C:{\cal G}(T) \to T$, the
evaluation function of $T$, where ${\cal G}(T)$ is the set of
contraction graphs of $T$.

A \emph{contraction graph} of $T$ is an oriented graph $G$ whose edges
and vertices are labeled in a certain way.  $G$ can be a graph in the
broadest sense of the word:  $G$ need not be connected (and may be
empty), an edge may go from a vertex to of vertices, an edge may go
from a vertex to itself, and there may be edges whose heads or tails
(or both) are not attached to any vertex.  $G$ may also have degenerate
``edges'' which are oriented circles without heads or tails at all:
$$\includegraphics{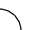}$$
The vertices of $G$ should be labeled by elements of $T$.  At a
vertex labeled by $t$, there should be $i(t)$ inward edges and $o(t)$
outward edges, the inward edges should be numbered from 1 to $i(t)$,
and similarly for the outward edges.  We also number the free heads of
edges of $G$ from 1 to $o$ for some $o$, and the free tails of edge of
$G$ from 1 to $i$ for some $i$.  We call $i$ the in-degree of $G$ and
$o$ the outdegree.  In short, a contraction graph looks very much like
a tensorial expression in arrow notation, where the ``tensors'' are
elements of $T$: 
$$\includegraphics{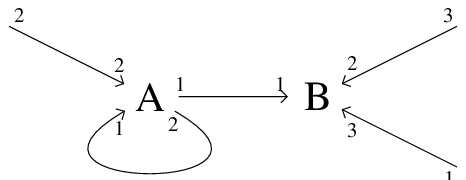}$$

The evaluation function $C$ should map a contraction graph $G$ to an
element of $T$ of the same in-degree and outdegree as $G$.  Moreover,
it should satisfy an axiom of substitution:  If $G$ and $H$ are
contraction graphs and $v$ is a vertex with the same in-degree
and outdegree as $H$, we may define the \emph{composition} $G \circ_v H$
by replacing the vertex $v$ by the graph $H$: $$\includegraphics{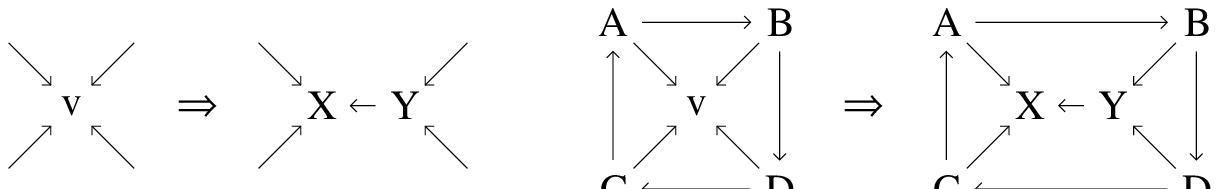}$$
The axiom of substitution says that if $C(H)$ is the label of $v$,
then $C(G \circ_v H) = C(G)$.

We give names to two distinguished tensors which exist in any
contraction category.  We call the value of the empty graph
the \emph{identity scalar} (or one), and the value of a circular
edge the \emph{dimension scalar}.

If $T$ and $U$ are contraction categories, a morphism from $T$ to $U$
is a function from the tensors of $T$ to the tensors of $U$ which respects
in-degree and outdegree and which commutes with the evaluation map.

The term ``contraction category'' may seem inappropriate since we
have not defined it as a kind of category.  However, we can view it
a category whose objects are ordered lists of left-pointing
arrows and right-pointing arrows (arrows in the sense of contraction 
graphs, that is; we call category-theoretic arrows ``morphisms''),
and whose morphisms are contraction graphs with an ordered list of
arrows on the left side and another ordered list on the right side:
$$\includegraphics{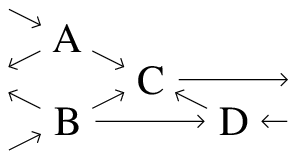}$$
We could modify the definition of a contraction category in natural
ways:  We could replace occurrences of the word ``set'' by ``class'',
and we could label the arrows by elements of a new set (or class) $A$.
In its full generality, a contraction category appears to be the same
as a strict, symmetric, compact, monoidal category \cite{JS:braided}.

As an example of a contraction category, we may define $T(V)$ to be the set of
all tensors over a finite-dimensional vector space $V$ with ground
field $k$.  In other words, $T = k \cup V \cup V^* \cup V \tensor V
\cup \ldots$.  The evaluation function for contraction graphs may be
constructed from the usual contraction and tensor product operations.
We see that the value of the dimension scalar is the dimension of $V$.
If $T$ is any contraction category, one can consider a
vector space $V$ together with a morphism $T \to T(V)$ to be a
representation of $T$.

A \emph{Hopf contraction category} is a contraction category with distinguished
elements called $M$, $\D$, $S$, $\eps$, and $i$ which satisfy the usual
identities of a Hopf algebra.  We may define the \emph{universal Hopf category}
$U$ to be the contraction category generated by these five letters with the
axioms they satisfy interpreted as relations.  That is, the elements
of $U$ are equivalence classes of contraction graphs on these five
letters, where two contraction graphs $G$ and $H$ are equivalent if
one can turn $G$ into $H$ by using the axioms of a Hopf algebra
as replacement rules.

We can equivalently say that a Hopf contraction category is a
contraction category $T$ together with a morphism from the universal
Hopf category to $T$.  A Hopf contraction category which happens to be
a vector space is a Hopf algebra.  We can also define the universal,
involutory Hopf category and the universal, involutory Hopf category of
invertible dimension in which trace is a left and right integral,
cotrace is a cointegral, and the dimension scalar has an inverse.
Anticipating the next section, we call the latter category the {\bf
Heegaard category}, written $\cal H$.  One can check that the Heegaard
category satisfies all of the identities of involutory Hopf algebras of
invertible dimension presented in section~\ref{Hopf}.

We can define the ``invariant'' $\sharp(M,{\cal H})$ as before, although its
value is not a number, but rather a formal expression which might be as
difficult to analyze as the original manifold.  We will see that
$\sharp(M,{\cal H})$ has the surprising property of being a complete
invariant on closed, irreducible 3-manifolds.

\section{$\sharp(M,{\cal H})$ as a complete formal invariant}

The value of $\sharp(M,{\cal H})$ involves factors of the dimension scalar for
normalization.  For simplicity we pass to the augmented Heegaard model
$\t{{\cal H}}$, in which we assume that the dimension scalar equals the
identity, and we consider $\sharp(M,\t{{\cal H}})$.

Section~\ref{invariant} gives a procedure for converting a generalized
Heegaard diagram $D$ to a contraction graph over the generating tensors
of $\t{{\cal H}}$.  We wish to construct an inverse to this procedure.

\begin{proposition} \label{inversep} Each element of $\t{{\cal H}}$ is
realized as $\sharp(M,\t{{\cal H}})$ for some bordism $M$.
\end{proposition}
\begin{proof}
Let $C$ be a graph composed of the generators of $\t{{\cal H}}$.  We
may re-write $i^a$, $e_a$, $M_{ab}{}^{c}$ and $\D^{ab}{}_c$ using
Corollary~\ref{tracegen}:
$$\includegraphics{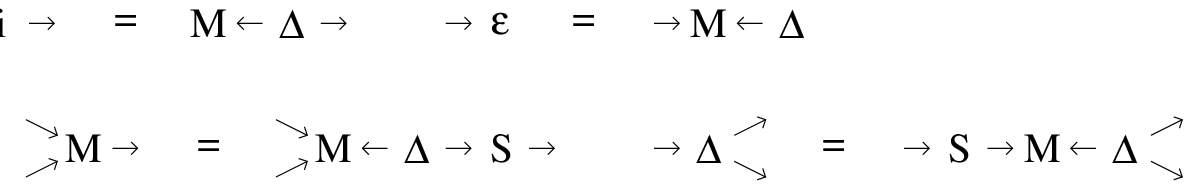}$$
Substituting these expressions into $C$, we get a new graph $C'$ composed of
the tracial product and coproduct tensors and the antipode.

\begin{figure}[ht]
\begin{center}
\subfigure[]{\includegraphics{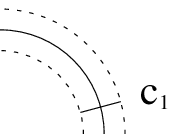}}
\subfigure[]{\includegraphics{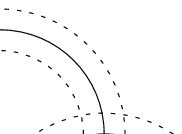}}
\caption{\label{f:core}}
\end{center}
\end{figure}

We build a Heegaard diagram using the graph $C'$.  For each tracial
product tensor in $C'$, we construct an oriented annulus with an
oriented core, which we call a ``lower circle'', together with a
labeled point for each index of the tensor, as in \fig{f:core}(a). We do the
same for each tracial coproduct, where we label the core as an ``upper
circle''.  If a product tensor and a coproduct tensor are contracted
together, the annuli should be identified, as in \fig{f:core}(b). The annuli
are identified in such a way that their orientations are compatible,
and so that at the crossing point the tangent vectors to the lower
circle and the upper circle, in that order, should be a positive
basis.  If a product tensor and a coproduct tensor are contracted
together with the antipode map in between, we perform the same gluing
operation, except we arrange the tangents to the lower and the upper
circle to form a negative basis.

The result of gluing these annuli together is a compact, oriented
surface $S$ which may have many components and many boundary
components.  To convert $S$ to a Heegaard diagram $D$ on a closed
surface, we glue it to a surface $T$ which has the same boundary as
$S$.  The result describes some bordism $M$.  By construction,
$\sharp(M,\t{{\cal H}}) = C'$.
 \end{proof}

There is an ambiguity in the construction of the Heegaard diagram $D$,
because we can choose many different surfaces $T$ with a given
boundary.  We may express this ambiguity by an equivalence relation.
Suppose that $D$ is a Heegaard diagram and $C$ is a circle on the
surface of $D$ which is disjoint from all of the upper and lower
circles of $D$.  Then we can cut the surface of $D$ along $C$ and
attach disks along the two resulting boundary components.  Call this a
\emph{circle move} on Heegaard diagrams.  Suppose instead that the
surface of $D$ has component which is a 2-sphere with no upper or lower
circles.  Then we can call the deletion of this component a \emph{sphere
move} on $D$.  Proposition~\ref{inversep} assigns a unique equivalence
class (under the sphere and circle moves) of Heegaard diagrams to each
contraction graph $C$.

By extension, we can say that two bordisms $M_1$ and $M_2$ differ by an
elementary move if there exist diagrams for $M_1$ and $M_2$ which
differ by a sphere or a circle move.  We can then consider equivalence
classes of bordisms.  It is easy to check that if $M$ is a bordism,
then the bordisms produced by Proposition~\ref{inversep} from
$\sharp(M,\t{{\cal H}})$ is equivalent to $M$.  One can also check that if
two contraction graphs $C_1$ and $C_2$ differ by one of the defining
relations of $\t{{\cal H}}$, then the corresponding bordisms are
equivalent (thereby establishing a converse to the proof of invariance
in section~\ref{invariant}).  Thus, $\sharp(M_1,\t{{\cal H}}) =
\sharp(M_2,\t{{\cal H}})$ if and only if $M_1$ and $M_2$ are equivalent
under sphere and circle moves.

\begin{figure}[ht]
\begin{center}
\includegraphics{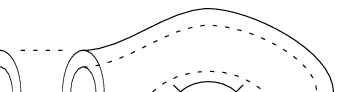}
\caption{\label{f:annulusmove}}
\end{center}
\end{figure}

We examine the effect of sphere and circle move on the topology of a
bordism $M$.  The sphere move corresponds to discarding a component of
$M$ which is a 3-sphere; we may call this move on manifolds a {\bf
3-sphere move}.  The circle move is more complicated.  Let $D$ be a
diagram for $M$ and let $S$ be the surface of $D$.  Recall that $M$ is
obtained by gluing disks and then balls to $S \cross [0,1]$.  Let $C$
be a circle on the surface of $D$ suitable for a circle move.  The
circle move corresponds to cutting $M$ along the annulus $C \cross
[0,1]$, and then gluing a 2-handle to each of the two annular scars
that result (see \fig{f:annulusmove}).  The operation may produce spherical
boundary components, which are irrelevant by the definition of a
bordism.  We call this operation an \emph{annulus move} on $M$.

It may happen that one or both of the circles $C \cross \{0\}$
and $C \cross \{1\}$ is trivial in its respective boundary component.
In this case the annulus move simplifies $M$ and we give it a 
second name as follows:

Suppose that only one of the circles $C \cross \{1\}$ and $C
\cross \{0\}$ is non-trivial, and assume that it is $C \cross \{1\}$.
We choose an embedded disk $S$ in the lower handlebody which bounds $C
\cross \{0\}$.  Then the circle move corresponds to compressing the
boundary of $M$ along the disk $S \cup C \cross [0,1]$.  We call this a
\emph{disk move} on $M$.

If both circles are trivial, then we can find a disk $S_0$ in the
lower boundary and a disk $S_1$ in the upper boundary so that $S_0 \cup
C \cross [0,1] \cup S_1$ is a 2-sphere, and the circle move on the
diagram of $M$ corresponds to cutting along this sphere, i.e. the
inverse connect-sum operation.  We call this a \emph{2-sphere move} on $M$.

If a 2-sphere in $M$ bounds a ball, then the only effect of the
corresponding 2-sphere move is to add a 3-sphere component to $M$,
which is reversed by the 3-sphere move.  We call such a 2-sphere move
\emph{degenerate}.  We will also refer to disk moves and 2-sphere moves as
degenerate annulus moves, and to 2-sphere moves as degenerate disk
moves.

We would like to determine when two bordisms are rendered equivalent
by a sequence of annulus moves, disk moves, and non-degenerate 2-sphere
moves and the inverse of these moves (and we apply the 3-sphere move
whenever possible).  We first observe that these three moves are
\emph{destructive} in the sense that they simplify a bordism, and their
inverses are constructive; only degenerate 2-sphere moves are neither
destructive nor constructive.  Specifically, we define the \emph{height} of a
bordism $M$ to be the ordinal $c\omega + s$, where $c$, the \emph{circle
number} of $\bd M$, is the maximum number of disjoint, non-trivial,
non-parallel circles in the boundary of $M$, and $s$, the \emph{Kneser number}
of $M$, is the maximum number of disjoint, non-trivial, non-parallel
2-spheres in $M$.  (Since $M$ is actually an equivalence class of
3-manifolds and some members of this class may have 2-sphere boundary,
we must stipulate that a 2-sphere which is parallel to the boundary of
a representative of $M$ is trivial.)  We may take $c = \dim H_1(\bd M)
- \dim H_0(\bd M)$; see \cite{Hempel:book} for a proof that $s$ is finite.
Clearly, destructive moves decrease the height of a bordism.  Thus, any
sequence of destructive moves, or \emph{path of destruction}, is finite.  The
main result of this section is a uniqueness theorem:

\begin{theorem} If $M$ is a bordism, then any two maximal paths of
destruction produce the same result.
\end{theorem}

Thus if no destructive moves can be performed on two different bordisms
$M_1$ and $M_2$ (e.g. if $M_1$ and $M_2$ are closed and prime), then
they are either equivalent or $\sharp(M_1,\t{{\cal H}}) \ne
\sharp(M_2,\t{{\cal H}})$.

This result is a special case of the Jaco-Shalen-Johannson decomposition
theorem, as presented in \cite{JS:seifert} and \cite{Johannson:homotopy}.  We
give a self-contained proof here:

\begin{figure}[ht]
\begin{center}
\subfigure[]{\includegraphics{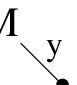}}
\subfigure[]{\includegraphics{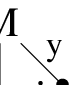}}
\caption{\label{f:paths}}
\end{center}
\end{figure}

\begin{proof} We outline the plan for a complicated proof by induction.  Let
$M$ be a counterexample of minimum height, and consider two possible
destructive moves $x$ and $y$ on $M$, and call the resulting manifolds $M_x$
and $M_y$.  If we can find two more moves $w$ and $z$ (which may or may not be
degenerate) such that $M_{xw} = M_{yz}$, then by induction any two maximal
paths of destruction beginning with $x$ and $y$ produce the same result, as
indicated in \fig{f:paths}(a).

The moves $x$ and $y$ involve cutting along surfaces $S_x$ and $S_y$. If these
two surfaces are disjoint, we can directly induct on the height of $M$, because
we can perform $x$ after $y$ and vice-versa, with $M_{xy} = M_{yx}$.  If $S_x$
and $S_y$ intersect, we put them in general position and induct on the number
of components of their intersection (each of which is an arc or a circle).  It
suffices to find another destructive move $z$ (possibly after switching $x$ and
$y$) such that $S_y \cap S_z$ has fewer components than $S_x \cap S_y$, and
there exist moves $v$ and $w$ such that $M_{zv} = M_{xw}$, as indicated in
\fig{f:paths}(b).

\begin{figure}[ht] \begin{center}
\subfigure[]{\includegraphics{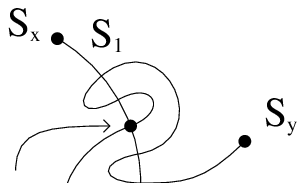}}
\subfigure[]{\includegraphics{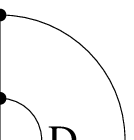}}
\subfigure[]{\includegraphics{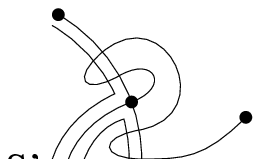}}
\caption{\label{f:innercircle}}
\end{center} \end{figure}

Suppose that one of the components of $S_x \cap S_y$ is a circle which
is contractible in one of $S_x$ and $S_y$, say $S_y$.  We may choose
an innermost circle $C \subset S_x \cap S_y$ inside this one. $C$
bounds a disk $D \subset S_y$ whose interior does not intersect $S_x$,
and $C$ divides $S_x$ into two surfaces $S_1$ and $S_2$, as shown in
\fig{f:innercircle}(a).  The move $x$ consists of cutting along $S_x$ and
gluing in either 2-handles or balls, and in both cases we can find a
disk $D_x$ in the region attached to $S_x$ which makes a sphere with
$D$, as shown in \fig{f:innercircle}(b).  We choose two surfaces $S'_1$ and
$S'_2$ which run parallel to $S_1 \cup D$ and $S_2 \cup D$ but which do
not intersect them or each other, as shown in \fig{f:innercircle}(c).
Cutting along $S'_1$ and cutting along $S'_2$ are valid moves, and
cutting along both is equivalent to performing $x$ and then cutting
along the sphere $D \cup D_x$.  Therefore at least one of $S'_1$ and
$S'_2$ must be a destructive move, and we may assume it is $S'_1$.  At
the same time, $S'_1$ has less intersection with $S_y$ than $S_x$
does.  We let $z$ be the move corresponding to $S'_1$.

\begin{figure}[ht] \begin{center}
\includegraphics{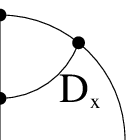}
\caption{\label{f:innerarc}}
\end{center} \end{figure}

From this point on, we suppose that there is no circular component of
$S_x \cap S_y$ which is contractible in either $S_x$ or $S_y$; in
particular, neither $x$ nor $y$ can be a 2-sphere move.  Suppose
instead that some component is an arc with both endpoints in upper or
lower boundary, say upper.  This arc, together with part of the
boundary of $S_y$, bounds a disk in $S_y$. We let $A \subset S_x \cap
S_y$ be the innermost arc in this disk.  $A$, together with part of the
boundary of $S_y$, bounds a disk $D \subset S_y$ whose interior does
not intersect $S_x$, and it divides $S_x$ into two surfaces $S_1$ and
$S_2$.  We argue as in the previous case, the only difference being
that $D \cup D_x$ is now a disk instead of a 2-sphere, as shown in
\fig{f:innerarc}.

\begin{figure}[ht] \begin{center}
\subfigure[]{\includegraphics{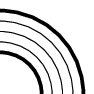}}
\subfigure[]{\includegraphics{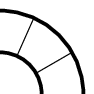}}
\caption{\label{f:annulus}}
\end{center} \end{figure}

\begin{figure}[ht] \begin{center}
\includegraphics{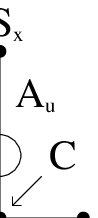}
\caption{\label{f:concentric}}
\end{center} \end{figure}

Thus, we may assume that both $x$ and $y$ are annulus moves.  There are
two possibilities:  Either $S_x$ and $S_y$ intersect in concentric
circles which lie between the upper and lower boundary of
$S_x$ and $S_y$, as in \fig{f:annulus}(a); or they intersect in radial
arcs, each arc connecting the upper and lower boundary of $S_x$ and
$S_y$, as in \fig{f:annulus}(b).  In the first case, let $C$ be the circle of
$S_x \cap S_y$ which is uppermost in $S_y$.  $C$ bounds an annulus
$A$ in $S_y$ which does not intersect $S_x$ in the interior.  It
divides $S_x$ into an upper annulus $A_u$ and a lower annulus $A_l$, as
in \fig{f:concentric}.  Cutting along $A \cup A_l$ is necessarily a
destructive move, because its upper boundary is a non-trivial circle.
We call this move $z$.  The circle $C$ bounds a disk $D_z$ in the
2-handle which is attached to $A \cup A_l$ in the move $z$, and $C$
bounds a disk $D_x$ in the 2-handle attached in the move $x$.  We see
that the move $z$ followed by compressing the disk $D_x \cup A$ yields
the same result as the move $x$ followed by compressing the disk $D_z
\cup A_l$.  As before, if we move $A \cup A_l$ away from the annulus
$A$, it intersects $S_y$ less than $S_x$ does.

Finally, suppose that $S_x$ and $S_y$ intersect in radial arcs.  Let
$\cal A$ be the boundary of a tubular neighborhood of $S_x \cup S_y$.
$\cal A$ is the union of disjoint annuli which do not intersect either
$S_x$ or $S_y$.  If we perform all of the corresponding annulus moves,
then $S_x$ and $S_y$ lie in a component of the resulting bordism which
is equivalent to $S \cross [0,1]$ for some closed surface $S$.  Any
maximal path of destruction removes this component entirely.  Thus,
$S_x$ and $S_y$ are rendered equivalent by the annulus moves
represented by $\cal A$.
 \end{proof}

\section{Unresolved questions}

Although $\sharp(M,{\cal H})$ is a complete formal invariant for closed
3-manifolds, it is not necessarily a complete invariant in the
computational sense, because it may be as difficult to distinguish
elements of ${\cal H}$ as it is to distinguish 3-manifolds by the
standard means.

\begin{conjecture} The Heegaard category ${\cal H}$ is residually
representable.  That is, if $a$ and $b$ are distinct elements of ${\cal
H}$, there exists a vector space $V$ and a morphism $m:{\cal H} \to
T(V)$ such that $m(a) \ne m(b)$.  \end{conjecture}

If this conjecture is true, $\sharp(M,{\cal H})$ is a complete invariant in
a more meaningful sense.  A few special cases of this conjecture might
be easy to settle.  Is there an involutory Hopf algebra $H$ for which
$\sharp(L(3,1),H) \ne \sharp(L(3,2),H)$, i.e. $\D^{abc}M_{abc} \ne
\D^{cba}M_{abc}$?  Even if $\sharp(M,H)$ is insensitive to orientation,
perhaps there is a Hopf algebra that distinguishes $L(7,1)$ from
$L(7,2)$, or $L(8,1)$ from $L(8,3)$ (correspondingly,
$\D^{abcdefg}M_{abcdefg}$ from $\D^{abcdefg}M_{acegbdf}$ and
$\D^{abcdefgh}M_{abcdefgh}$ from $\D^{abcdefgh}M_{adgbehcf}$).

The construction given in this paper is almost certainly not as general as it
could be.  For example, Dijkgraaf \cite{Dijkgraaf:thesis} mentions the
following state model on triangulations which was described to me by Vaughan
Jones: As before, we choose a finite group $G$ to be the state set, we assign
states to the oriented edges of a triangulation, and we assign an interaction
to each face which is 0 unless the product of the states of the face's edges is
the identity.  However, we also have an interaction for each tetrahedron, which
we assume is non-zero whenever the value of the tetrahedron's faces is
non-zero. We may view this interaction as a complex-valued function
$I(a,b,c,d,e,f)$, with $a,b,c,d,e,f \in G$, and therefore as a cochain on the
canonical triangulation of $K(G,1)$ with coefficients in the multiplicative
group $\C^*$.  The constraint that the state model is a topological invariant
is equivalent to the condition that this cochain is actually a cocycle, and two
different cocycles yield the same invariant if they are cohomologous.  Thus, we
get an invariant for every pair consisting of a finite group and an element of
$H^3(G,\C^*)$.  If we choose the trivial cohomology class, we get
$\sharp(M,\C[G])$.  Is there a mutual generalization of this invariant and
$\sharp(M,H)$?

Can the invariant be extended to non-involutory Hopf algebras?  As was
mentioned in section~\ref{prop}, $\sharp(M,H)$ corresponds to a certain link
invariant involving the quantum double of $H$.  Can it be generalized
to quasi-triangular Hopf algebras other than the quantum double?


\begin{thebibliography}{10}

\bibitem{Baxter:exactly}
Rodney~J. Baxter, \emph{Exactly solved models in statistical mechanics},
  Academic Press, London, 1982.

\bibitem{Bergman:personal}
George Bergman, personal communication.

\bibitem{Cerf:strat}
Jean Cerf, \emph{La stratification naturelle des espaces de fonctions
  diff\'erentiables r\'eelles et le th\'eor\`eme de la pseudo-isotopie}, Inst.
  Hautes \'Etudes Sci. Publ. Math. (1970), no.~39, 5--173.

\bibitem{Dijkgraaf:thesis}
Robbert Dijkgraaf, \emph{A geometrical approach to two dimensional conformal
  field theory}, Ph.d. thesis, Rijksuniversiteit Utrecht, 1989.

\bibitem{Drinfeld:quantum}
Vladimir Drinfel'd, \emph{Quantum groups}, Proceedings of the International
  Congress of Mathematicians (Berkeley, California), vol.~1, 1986,
  pp.~798--820.

\bibitem{HOMFLY:invariant}
Peter Freyd, David Yetter, Jim Hoste, W.~B.~R. Lickorish, Kenneth Millett, and
  Adrian Ocneanu, \emph{A new polynomial invariant of knots and links}, Bull.
  Amer. Math. Soc. (N.S.) \textbf{12} (1985), no.~2, 183--312.

\bibitem{Hempel:book}
John Hempel, \emph{$3$-{Manifolds}}, Ann. of Math. Studies, vol.~86, Princeton
  University Press, Princeton, NJ, 1976.

\bibitem{JS:seifert}
William~H. Jaco and Peter~B. Shalen, \emph{Seifert-fibered spaces in
  3-manifolds}, Mem. Amer. Math. Soc. \textbf{21} (1979), no.~220, viii+192.

\bibitem{Johannson:homotopy}
Klaus Johannson, \emph{Homotopy equivalences of $3$-manifolds with boundaries},
  Springer, Berlin, 1979.

\bibitem{Jones:pacific}
Vaughan F.~R. Jones, \emph{On knot invariants related to some statistical
  mechanical models}, Pacific J. Math. \textbf{137} (1989), no.~2, 311--334.

\bibitem{JS:braided}
Andr{\'e} Joyal and Ross Street, \emph{Braided tensor categories}, Adv. Math.
  \textbf{102} (1993), no.~1, 20--78.

\bibitem{Kauffman:state}
Louis~H. Kauffman, \emph{State models and the {Jones} polynomial}, Topology
  \textbf{26} (1987), no.~3, 395--407.

\bibitem{KR:links}
Anatol~N. Kirillov and Nicolai~Yu. Reshetikhin, \emph{Representations of the
  algebra ${U_q(sl(2))}$, $q$-orthogonal polynomials, and invariants of links},
  Infinite-dimensional {Lie} algebras and groups ({Luminy-Marseille}, 1988),
  World Sci. Publishing, 1989, pp.~285--339.

\bibitem{LR:semisimple0}
Richard~G. Larson and David~E. Radford, \emph{Finite-dimensional cosemisimple
  {Hopf} algebras in characteristic 0 are semisimple}, J. Algebra \textbf{117}
  (1988), no.~2, 267--289.

\bibitem{RT:manifolds}
Nicolai~Yu. Reshetikhin and Vladimir~G. Turaev, \emph{Invariants of 3-manifolds
  via link polynomials and quantum groups}, Invent. Math. \textbf{103} (1991),
  no.~3, 547--597.

\bibitem{Sweedler:hopf}
Moss~E. Sweedler, \emph{Hopf algebras}, W. A. Benjamin, Inc., New York, 1969.

\bibitem{Wald:general}
Robert~M. Wald, \emph{General relativity}, University of Chicago Press,
  Chicago, IL, 1984.

\end{thebibliography}

\providecommand{\bysame}{\leavevmode\hbox to3em{\hrulefill}\thinspace}

\end{document}